\documentclass[10pt,reqno,final]{amsart}
\usepackage{amsfonts}
\usepackage[notcite,notref]{showkeys}
\usepackage{url,hyperref,multirow}
\usepackage{color}
\usepackage{cases}
\usepackage{subfigure}
\usepackage{epsfig,amsmath,bm,caption2,epstopdf}
\usepackage{amssymb,version,graphicx,fancybox,mathrsfs,pifont,booktabs}

\catcode`\@=11 \theoremstyle{plain}

\@addtoreset{figure}{section}
\renewcommand\thefigure{\@arabic\c@figure}
\@addtoreset{table}{section}
\renewcommand\thetable{\@arabic\c@table}

\newtheorem{thm}{\bf Theorem}
\newtheorem{cor}{\bf Corollary}

\newtheorem{lmm}{\bf Lemma}

\theoremstyle{remark}
\newtheorem{remark}{\bf Remark}

\def \ri {{\rm i}}
\def \rd {{\rm d}}

\newcommand{\bs}[1]{\boldsymbol{#1}}

\usepackage{graphicx}

\begin{document}

%
%
%
%

\title [FEM for the integral fractional Laplacian]
{On explicit form of the FEM stiffness matrix for the integral fractional Laplacian on non-uniform meshes}

\author[H. Chen,\, C. Sheng\, and\,    L.-L.  Wang  ]{Hongbin Chen${}^{1}$,\;\; Changtao Sheng${}^{2}$ \;and\; Li-Lian Wang${}^{2}$  }

\thanks{${}^{1}$Institute of Mathematics and Physics, College of Science, Central South University of Forestry and Technology, Changsha, Hunan 410081, China.  Email: hongbinchen@csuft.edu.cn (H. Chen).\\
\indent ${}^{2}$Division of Mathematical Sciences, School of Physical and Mathematical Sciences, Nanyang Technological University (NTU), 637371, Singapore. The research of the authors is partially supported by Singapore MOE AcRF Tier 2 Grants:  MOE2018-T2-1-059 and MOE2017-T2-2-144. Emails: ctsheng@ntu.edu.sg (C. Sheng) and  lilian@ntu.edu.sg (L. Wang).
\\ \indent The first author would like to thank NTU for hosting his visit devoted to  this collaborative work.  
         }

 \begin{abstract} We derive  exact form of the piecewise-linear finite element stiffness matrix on general non-uniform meshes for  the integral fractional Laplacian operator in one dimension, where the derivation is 
 accomplished in the Fourier transformed space. With such an exact formulation at our disposal, 
 we are able to  numerically study some intrinsic properties of the fractional stiffness matrix on some commonly used non-uniform meshes (e.g., the graded mesh), in particular, to examine their seamless transition to those of the usual Laplacian.  
\end{abstract}



\maketitle

\vspace*{-18pt} 

\section{Introduction}

There has been a burgeoning of recent interest in nonlocal and fractional models, largely due to  the advancement in 
both computing power and  computational algorithms.  
The integral fractional Laplacian (IFL) is deemed as one of the most prominent nonlocal operators, but unfortunately, it poses 
 more  challenges  in numerical solutions of the related models. 
 Among very  limited works on finite element approximation of the IFL,
 D'Elia and Gunzburger \cite{d2013fractional}  considered the FEM discretisation on  non-uniform meshes in one dimension.
 The entries of the FEM stiffness matrix  therein were computed by the Gauss quadrature rule, and the adaptive Gauss–Kronrod quadrature (with a built-in function in Matlab) was resorted to  approximate  the double integrals with singular kernels when the mesh size is small.  
 
 In this paper, we compute the entries of the stiffness matrix  in the 
  Fourier transformed space  based on the definition of the IFL:  for $s\ge 0,$ 
  \begin{equation}\label{Ftransform}
\begin{split}
&(-\Delta)^s u(x):={\mathscr F}^{-1}\big[|\xi|^{2s}  {\mathscr F}[u](\xi)\big](x),
\end{split}
\end{equation}
where $u(x)$ on  $\mathbb R=(-\infty, \infty)$ is of Schwartz class, and $\mathscr{F}$ denotes the Fourier transform with the inverse ${\mathscr F}^{-1}.$  In fact,  for $s\in (0,1),$ the IFL of $u(x)$ can be equivalently defined by 
 \begin{equation}\label{fracLap-defn}
(-\Delta)^s u(x)=C_s\, {\rm p.v.}\! \int_{\mathbb R} \frac{u(x)-u(y)}{|x-y|^{1+2s}} \rd y\;\;\;\;  
 \text{with} \quad C_s=\frac{2^{2s}s\,\Gamma(s+\frac{1}{2})}{\sqrt \pi \, \Gamma(1-s)},
\end{equation}
where ``p.v." stands for the principle value and $\Gamma(\cdot)$ in the Gamma function.
  As opposed to  \cite{d2013fractional} and limited existing works implemented via \eqref{fracLap-defn}, 
  the use of the formulation \eqref{Ftransform} enables us to evaluate  the entries explicitly. With such an analytic representation, we can study  some intrinsic properties of the stiffness matrix and related numerical issues when the meshes are highly non-uniform. 

\section{Main result}

Consider the model equation with  a global homogeneous Dirichlet boundary condition:
\begin{equation}\label{dDsLap00}
(-\Delta)^{s} u(x)=f(x), \quad x\in\Omega=(a,b); \quad u(x)=0, \quad x\in\Omega^c=\mathbb{R}\backslash\bar{\Omega},
\end{equation}
where $f(x)$ is a given continous function. 
Let $\{\phi_j(x)\}_{j=1}^{N-1}$ be the piecewise linear finite element basis (i.e., the standard ``hat'' functions) associate with the  partition  
\begin{equation}\label{partition}
\Omega:\;\;a=x_0<x_1<\cdots<x_{N}=b, 
\end{equation}
and satisfy $\phi_j(x)\equiv0$ for $x\in\Omega^c$.  The piecewise linear FEM approximation to \eqref{dDsLap00} is to find $u_h\in V_h^0={\rm span}\{\phi_j\,:\, 1\le j\le N-1\}$ such that  
\begin{equation}\label{dDsLap02}
a_s(u_h, v_h)=\int_{\mathbb R} ((-\Delta)^{s/2} u_h) ((-\Delta)^{s/2} v_h) {\rm d} x =\int_{\Omega} f(x) v_h(x)\, {\rm d} x, \quad \forall\, v_h\in V_h^0,
\end{equation}
which admits  a unique solution by the standard Lax-Milgram lemma.  

\subsection{Main result}
Our main purpose is to show that the stiffness matrix, denoted by $\bs S,$ with the entries 
\begin{equation*}\label{Sjkform}
S_{kj}=S_{jk}= a_s(\phi_j, \phi_k),  
\quad 1\le k,j\le N-1,
\end{equation*}
has  the following explicit form.  
\begin{thm}\label{thmmain}  If $s\in(0,3/2)$ but $s\not=1/2,$ then the entries of the stiffness matrix $\bs S$ 
  can be explicitly evaluated by 
\begin{equation}\label{HY2.3}
S_{jk} = \widehat C_s\, {\bs c_j\, \bs{D}^{k}_{j}\,\bs c_k^t}\;\;\;\;  {\rm with}\;\;\; \widehat C_s:= \frac 1 {2\Gamma(4-2s)\cos(s\pi)},
\end{equation}
where  
\begin{equation}\label{dkj}\begin{split}
&\bs c_l=\Big(\frac{1}{h_l}, -\frac{1}{h_l}-\frac{1}{h_{l+1}}, \frac{1}{h_{l+1}}\Big),\quad \bs{D}^{k}_{j}=\!\begin{pmatrix}
 (d^{k-1}_{j-1})^{\gamma}&(d^k_{j-1})^{\gamma} & (d^{k+1}_{j-1})^{\gamma}\\[6pt]
( d^{k-1}_{j})^{\gamma} &(d^k_{j})^{\gamma} &(d^{k+1}_{j})^{\gamma}\\[6pt]
 (d^{k-1}_{j+1})^{\gamma}&(d^k_{j+1})^{\gamma} & (d^{k+1}_{j+1})^{\gamma}
  \end{pmatrix},
\end{split}\end{equation}
with $h_\ell=x_{\ell}-x_{\ell-1},$ $ d_\iota^\ell=|x_\iota-x_\ell|,$ and $\gamma=3-2s.$

If $s=1/2$,  then  we have 
\begin{equation}\label{StiffMatrix01/2}
  S_{jk}= \frac{1}{2\pi} \bs c_j  \begin{pmatrix}
 (d^{k-1}_{j-1})^2\ln d^{k-1}_{j-1}&(d^k_{j-1})^2\ln d^k_{j-1} & (d^{k+1}_{j-1})^2\ln d^{k+1}_{j-1}\\[6pt]
( d^{k-1}_{j})^2\ln d^{k-1}_{j}&(d^k_{j})^2\ln d^k_{j} &(d^{k+1}_{j})^2\ln d^{k+1}_{j}\\[6pt]
 (d^{k-1}_{j+1})^2\ln d^{k-1}_{j+1}&(d^k_{j+1})^2 \ln d^k_{j+1}& (d^{k+1}_{j+1})^2\ln d^{k+1}_{j+1}
  \end{pmatrix}\bs c_k^t,
\end{equation}
where we understand  $(d^\ell_\iota)^2\ln d^\ell_\iota=0$ when $d^\ell_\iota=0.$
\end{thm}

Prior to the proof, we discuss some implications and consequences of the main result.  Observe from the above that for fixed $s\in (0,3/2),$ the matrix $\bs S$ is completely determined by the partition  \eqref{partition}, and  for fixed $j,k$,  the entry $S_{jk}$  only involves the  grid points: $\{x_{j+p}\}_{p=0,\pm 1}$ and $\{x_{k+q}\}_{q=0,\pm 1}.$ Interestingly,  
$S_{jk}$ turns out be a finite difference approximation of 
\begin{equation*}\label{hatdxy}
\hat d(x,y)=\widehat C_s\, |x-y|^{3-2s},\;\; {\rm if}\;\;  s\not=\frac 1 2;\quad  \hat d(x,y)=\frac 1 {2\pi}(x-y)^2\ln|x-y|,\;\; {\rm if}\;\;  s=\frac 1 2.
\end{equation*}
Indeed, one verifies from Theorem \ref{thmmain} the following alternative representation. 
\begin{cor}\label{FDmax}   For $s\in (0,3/2),$ the entry $S_{jk}$ can be written as a finite difference form
\begin{equation}\label{newform}
S_{jk}=\delta_y^2\delta_x^2  \hat d_j^k=\delta_x^2\delta_y^2  \hat d_j^k,\quad 1\le j,k\le N-1,
\end{equation}
where  $\hat d_j^k=\hat d(x_j,x_k)$ and 
\begin{equation*}\label{newform2}
\delta_x^2 \hat d_{j}^k:= \frac{\hat d_{j}^k-\hat d_{j-1}^k}{h_j}- \frac{\hat d_{j+1}^k-\hat d_{j}^k}{h_{j+1}},\quad 
\delta_y^2 \hat d_{j}^k:= \frac{\hat d_{j}^k-\hat d_{j}^{k-1}}{h_k}- \frac{\hat d_{j}^{k+1}-\hat d_{j}^k}{h_{k+1}}.
\end{equation*}
\end{cor}
\begin{remark}
 It is noteworthy that from standard finite difference formula, we have 
 \begin{equation}\label{fd2form}
 \delta_x^2 \hat d_{j}^k=\frac1 2(h_j+h_{j+1}) \partial_x^2 \hat d(x_j,x_k)+O(h_j^2+h_{j+1}^2).
 \end{equation}
Thus, $S_{jk}$ is a nine-point finite difference approximate of $\hat d(x,y)$ on  $\{(x_{j+p}, x_{k+q})\}_{p,q=0,\pm 1}.$  \qed
\end{remark}

\begin{remark}
When $s\to 0$ and $s=1$, the matrix $\bs S$ in  Theorem {\rm \ref{thmmain}}  reduces to  the usual (tridiagonal) FEM mass matrix 
$\bs M={\rm diag}(h_j/6, (h_j+h_{j+1})/3, h_{j+1}/6)$ and the  stiffness matrix $\bs S={\rm diag}(-1/h_j, 1/h_j+1/h_{j+1},-1/h_{j+1})$, 
respectively. If \eqref{partition} is a uniform partition of $(a,b)$ with $h=h_j$, then \eqref{HY2.3}  reduces to  
 \begin{equation*} 
 S_{jk}=\widehat C_s {h^{1-2s}}  \sum_{i=-2}^2 w_i \big||k-j|+i\big|^{3-2s} \;\; {\rm with}\;\;   w_0=6,\;\; w_{\pm 1}=-4,\;\; w_{\pm 2}=1.
\end{equation*}
In this case, the stiffness matrix $\bs S$ is a Toeplitz matrix  (cf.  \cite{liu2020diagonal} and also for some other interesting properties).
\qed
\end{remark}
%

\subsection{Proof of Theorem \ref{thmmain}}
Recap on the piecewise linear FEM basis associated with \eqref{partition}:  
\begin{equation}\label{NonuniformP}
\phi_{j} (x) =
\begin{cases}
 c_j (x-x_{j-1}),\quad & x\in (x_{j-1},x_{j}),\\
 c_{j+1}(x_{j+1} - x),\quad & x\in (x_{j},x_{j+1}),\\
 0, \quad & \text{elsewhere on} \,\, \mathbb{R},
\end{cases}\quad  c_\ell=\frac1 {h_\ell}, 
\end{equation}
for $1\le j\le N-1.$ Using   integration by parts leads to 
\begin{equation}\label{HY2.2}
\begin{aligned}
\mathscr{F}[\phi_{j}](\xi)
& = \frac{1}{\sqrt{2\pi}}\int_{\mathbb R} \phi_{j}(x)e^{-\ri x\xi}\,{\rm d}x
  = \frac{1}{\sqrt{2\pi}}\int_{x_{j-1}}^{x_{j+1}}\phi_{j}(x)e^{-\ri x\xi}\,{\rm d}x \\
& = \frac{ c_{j}}{\sqrt{2\pi}}\int_{x_{j-1}}^{x_{j}}(x-x_{j-1})e^{-\ri x\xi}\,{\rm d}x
+ \frac{c_{j+1}}{\sqrt{2\pi}}\int_{x_{j}}^{x_{j+1}}(x_{j+1}-x)e^{-\ri x\xi}\,{\rm d}x \\
& = \frac{1}{\sqrt{2\pi}}\Big{[}\frac{c_{j}}{\xi^{2}}(e^{-\ri x_{j}\xi}-e^{-\ri x_{j-1}\xi}) - \frac{c_{j+1}}{\xi^{2}}(e^{-\ri x_{j+1}\xi}-e^{-\ri x_{j}\xi})\Big{]}\\
& = -\frac{1}{\sqrt{2\pi}}\frac{ c_{j}e^{-\ri x_{j-1}\xi}-(c_{j}+ c_{j+1})e^{-\ri x_{j}\xi}+ c_{j+1}e^{-\ri x_{j+1}\xi}}{\xi^{2}}= -\frac{1}{\sqrt{2\pi} \xi^2} \bs c_j \bs e_j(\xi),
\end{aligned}
\end{equation}
 where $\bs e_j(\xi):=(e^{-\ri x_{j-1}\xi}, e^{-\ri x_{j}\xi}, e^{-\ri x_{j+1}\xi})^t.$
In view of \eqref{Ftransform}, \eqref{dDsLap02}, and \eqref{HY2.2}, we obtain from direct calculation and the parity of cosines and sines that
\begin{equation}\label{HY2.5}\begin{split}
 S_{jk}  &=  \int_{\mathbb{R}}|\xi|^{2s}\mathscr{F}[ \phi_{j}](\xi)\overline{\mathscr{F}[ \phi_{k}](\xi)}\,{\rm d}\xi 
= \frac{1}{2\pi} \bs c_j \Big(\int_{\mathbb R} |\xi|^{2s-4} \bs e_j(\xi) \bs e_k^t(-\xi) {\rm d}\xi \Big)\bs c_k^t   \\
 & 
 = \frac{1}{\pi}\int_{0}^{\infty}\xi^{2s-4}f_{jk}(\xi)\,{\rm d}\xi,
\end{split}
\end{equation}
where 
\begin{equation*}\label{nonuni1}
 f_{jk}(\xi) = \bs c_j \bs{F}_{jk}(\xi) \bs c_k^t, \quad 
 \bs{F}_{jk}(\xi)=\begin{pmatrix}
 \cos(d_{j-1}^{k-1}\xi) &  \cos(d_{j-1}^k\xi)&  \cos(d_{j-1}^{k+1}\xi)\\[2pt]
 \cos(d_{j}^{k-1}\xi) &  \cos(d_{j}^k\xi)&  \cos(d_{j}^{k+1}\xi)\\[2pt]
 \cos(d_{j+1}^{k-1}\xi) &  \cos(d_{j+1}^k\xi) &  \cos(d_{j+1}^{k+1}\xi)
  \end{pmatrix}.
\end{equation*}
One  verifies from direct calculation or   finite difference approximation  (cf.\,\eqref{newform} and \eqref{fd2form})   
that $f_{jk}(0)=f'_{jk}(0)=f''_{jk}(0) = 0.$
\smallskip 

We first consider  $s\in(1,3/2)$.  Recall the integral identity (cf. \cite[p. 440]{Gradshteyn2015Book}):
\begin{equation}\label{sinint}
\int^{\infty}_0x^{\mu-1}\sin(ax)\, {\rm d}x=\frac{\Gamma(\mu)}{a^\mu}\sin\Big(\frac{\mu\pi}{2}\Big),\quad a>0, \;\; \mu\in(0,1).
\end{equation}
We derive from \eqref{HY2.5} and integration by parts that
\begin{equation}
\begin{aligned}\label{nonuniform3}
\int_{0}^{\infty} & \xi^{2s-4}f_{jk}(\xi)\,{\rm d}\xi
 = \frac{1}{2s-3}\Big{\{}\xi^{2s-3}f_{jk}(\xi)\big|_{0}^{\infty} - \int_{0}^{\infty}\xi^{2s-3}f'_{jk}(\xi)\,{\rm d}\xi\Big{\}}\\
& = -\frac{1}{2s-3} \int_{0}^{\infty}\xi^{2s-3}f'_{jk}(\xi)\,{\rm d}\xi=  -\frac{1}{2s-3} \bs c_j 
\Big(\int_{0}^{\infty}\xi^{2s-3}\bs F'_{jk}(\xi){\rm d}\xi\Big)  \bs c_k^t.
\end{aligned}
\end{equation}
Applying \eqref{sinint} with $\mu=2s-2$  to each entry of  $\xi^{2s-3} \bs F'_{jk}(\xi)$ yields 
\begin{equation*}\label{HY2.10}
\int_{0}^{\infty}\xi^{2s-4}f_{jk}(\xi)\,{\rm d}\xi = -\Gamma(2s-3)\sin(s\pi)\,{\bs c_j  \bs{D}^{k}_{j} \bs c_k^t}.
\end{equation*}

 We next consider $s\in(1/2,1).$ 
 Recall that (cf. \cite[p. 441]{Gradshteyn2015Book})
\begin{equation}\label{cosint}
\int_0^{\infty} x^{\mu-1}\cos(ax){\rm d}x=\frac{\Gamma(\mu)}{a^\mu}\cos\Big(\frac{\mu\pi}{2}\Big),\quad a>0,\;  \mu\in(0,1).
\end{equation}
Applying integration by parts one more time to \eqref{nonuniform3}, we derive from \eqref{cosint}  with $\mu=2s-1$  that 
\begin{equation}
\begin{aligned}\label{nonuniform2}
\int_{0}^{\infty} \xi^{2s-4}f_{jk}(\xi)\,{\rm d}\xi
 & = \frac{1}{(2s-3)(2s-2)} \int_{0}^{\infty}\xi^{2s-2}f''_{jk}(\xi)\,{\rm d}\xi\\
&= \frac{1}{(2s-3)(2s-2)} \bs c_j 
\Big(\int_{0}^{\infty}\xi^{2s-2}\bs F''_{jk}(\xi){\rm d}\xi\Big)  \bs c_k^t\\
&= -\Gamma(2s-3)\sin(s\pi)\,{\bs c_j \bs{D}^{k}_{j}\bs c_k^t}.
\end{aligned}
\end{equation}
%

We now turn to  $s\in(0,1/2)$.  
Similarly, we integrate \eqref{nonuniform2} by parts once more
and use   \eqref{sinint} with $\mu=2s$ to obtain 
\begin{equation*}
\begin{aligned}\label{nonuniform1}
\int_{0}^{\infty} \xi^{2s-4}f_{jk}(\xi){\rm d}\xi
  &= -\frac{1}{(2s-3)(2s-2)(2s-1)} \int_{0}^{\infty}\xi^{2s-1}f'''_{jk}(\xi){\rm d}\xi\\
 & = -\Gamma(2s-3)\sin(s\pi) \,{\bs c_j \bs{D}^{k}_{j}\bs c_k^t}.
\end{aligned}
\end{equation*}
Then, using the property: $\Gamma(z)\Gamma(1-z)= \pi/{\sin \pi z}$ ($z\not=0,-1,\cdots$), we can reformulate the constant and then obtain the desired representation in 
\eqref{HY2.3} for all three cases. 

Finally, for $s=1/2,$  using the fact and the basic limit 
\begin{equation*}
 \lim_{s \to \frac{1}{2}} \textbf{c}_j \bs D_j^k \textbf{c}_k^t = 0, \quad \ln z = \lim_{\delta\to 0}\frac{z^{\delta}-1}{\delta}, \;\;\;  z>0,
\end{equation*}
%
we can directly take limit on \eqref{HY2.3}: 
$$
S_{jk}=\frac 1 4   \lim_{s\to \frac 12}\frac{\bs c_j \bs D_j^k \bs c_k^t}{\cos (s\pi)},
$$
and use the L'Hospital's rule to obtain \eqref{StiffMatrix01/2}. This completes the proof.

\section{FEM on graded meshes} 

It is known that the graded meshes are commonly used in finite element approximation of solutions with boundary singularities. 
In general, the mesh geometry affects not only the approximation error of the finite element solution but also the spectral properties of the corresponding stiffness matrix. It is a well-studied topic in the integer-order case, but much less known in this fractional setting.  

\subsection{A singular mapping} We propose to generate the graded mesh on $[a,b]$ for the solutions with singularities at the endpoint(s)  by the singular mapping \cite{wang2005error}: 
\begin{equation}\label{GradedMap}
x=g(y; \alpha, \beta) = a+(b-a) \frac{B(y; \alpha, \beta)}{B(\alpha, \beta)}\;\;\;  {\rm with}\;\;\; 
B(y; \alpha, \beta) = \int_0^yt^{\alpha-1}(1-t)^{\beta-1} {\rm d}t,
\end{equation}
for $y\in [0,1], x\in [a, b],$ and $\alpha,\beta\ge 1,$ where $B(y, \alpha, \beta)$ is  incomplete Beta function and $B(\alpha,\beta)=B(1; \alpha,\beta)$ is the Beta function. It is a one-to-one mapping such that $a=g(0; \alpha, \beta)$ and $b=g(1;\alpha,\beta).$ If $\alpha=\beta=1,$ it reduces to a linear transform.  Let $\{y_j=j/N\}_{j=0}^N$ be a uniform partition of the reference interval $[0,1].$ Then the mapped grids on $[a,b]$ are given by 
\begin{equation}\label{xgrids}
x_j:=x_{N,j}^{(\alpha,\beta)}= g(y_j; \alpha, \beta)=g(j/N; \alpha, \beta),\quad 0\le j\le N.
\end{equation}
By the mean value theorem, 
\begin{equation*}\label{xgrids2}
h_j=x_j-x_{j-1}=\frac{dx}{dy}\Big|_{y=\xi_j}(y_j-y_{j-1}) = \frac{b-a}{B(\alpha,\beta)}\frac{\xi_j^{\alpha-1}(1-\xi_j)^{\beta-1}}{N},
\end{equation*}
for some $\xi_j\in (y_{j-1}, y_j), 1\le j\le N.$ This implies 
$$
h_1\le \frac{b-a}{B(\alpha,\beta)} \frac 1 {N^\alpha},\quad h_N\le \frac{b-a}{B(\alpha,\beta)} \frac 1 {N^\beta},
$$ 
and the grid spacing near $x=a$ (resp. $x=b$) is of order $O(N^{-\alpha})$ (resp. $O(N^{-\beta})),$ while it remains $O(N^{-1})$  slightly away from the endpoints. 
\begin{remark}  If $\alpha>1$ and $\beta=1,$ then \eqref{xgrids} reduces to 
\begin{equation}\label{xjf}
x_j=g(y_j;\alpha,1)= a+(b-a) (y_j)^\alpha= a+(b-a) \Big(\frac {j} N\Big)^\alpha,\quad 0\le j\le N,
\end{equation}
which leads to a graded mesh with grid clustering near the left endpoint $x=a.$  Likewise, $\{x_j=g(y_j;1,\beta)\}$ with $\beta>1$ produces a graded mesh for the right end-point singularity.   It is noteworthy that the distribution of the mapped grids with $\alpha=\beta>1$
for symmetric end-point singularities  is slightly different from that generated by \eqref{xjf} and used in practice:
\begin{equation}\label{GradedM3}
x_j =\tilde g(y_j;\alpha)=
\begin{cases}
 a + \dfrac{b-a}{2}\Big(\dfrac{2j}{N}\Big)^{\alpha},  &  j = 0, 1, \cdots, N/2-1,\\[6pt]
 b - \dfrac{b-a}{2}\Big(2-\dfrac{2j}{N}\Big)^{\alpha}, \;\;  &  j = N/2, N/2+1, \cdots, N,
\end{cases}
\end{equation}
where $N>1$ is assumed to be an even integer, and  the underlying mapping  has a limited regularity at $x=(b+a)/2.$ 
However, this is not the case, if one uses \eqref{GradedMap}.
\qed
\end{remark}

 \subsection{Conditioning of the stiffness matrix}
  According to \cite[(26)]{fried1972condition}, the condition number of stiffness matrix $\bs S$ for an integer-order elliptic problem of the $2m$th order ($m=1$ harmonic, $m=2$ biharmonic) is given by
\begin{equation}\label{condm}
{\rm Cond}(\bs S)=c\big(h_{\max}/h_{\min}\big)^{2m-1}N^{2m},
\end{equation}
where $c$ is a positive numerical constant, $N$ is the degree of freedom, and  $h_{\max}, h_{\min}$ are the largest and smallest mesh sizes, respectively. It indicates a clear dependence of the condition number on  the mesh ratio $\rho:=h_{\max}/h_{\min},$ and 
the condition number is  greatly magnified for a highly non-uniform mesh, compared with  a quasi-uniform mesh with constant $\rho$. 

\begin{figure}[!th]
\subfigure[${\rm Cond}\sim N^{\alpha(2s-1)+1}$\hspace*{22pt}]{\hspace{-26pt}
\begin{minipage}[t]{0.40\textwidth}
\centering 
\rotatebox[origin=cc]{-0}{\includegraphics[width=0.8\textwidth,height=0.65\textwidth]{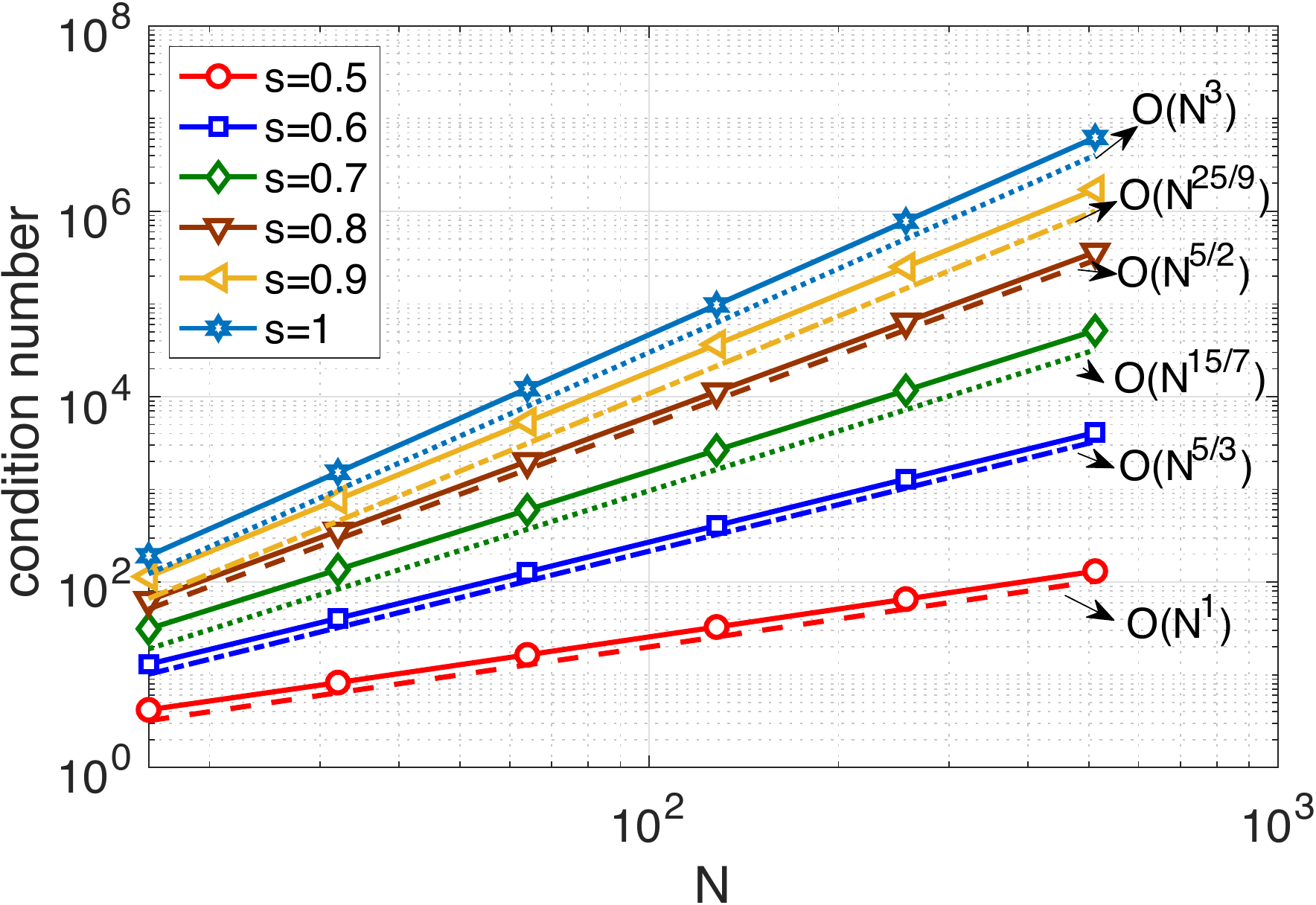}}
\end{minipage}}
\subfigure[$\lambda_{\min}\sim N^{-1}$\hspace*{24pt}]{\hspace{-28pt}
\begin{minipage}[t]{0.40\textwidth}
\centering 
\rotatebox[origin=cc]{-0}{\includegraphics[width=0.8\textwidth,height=0.65\textwidth]{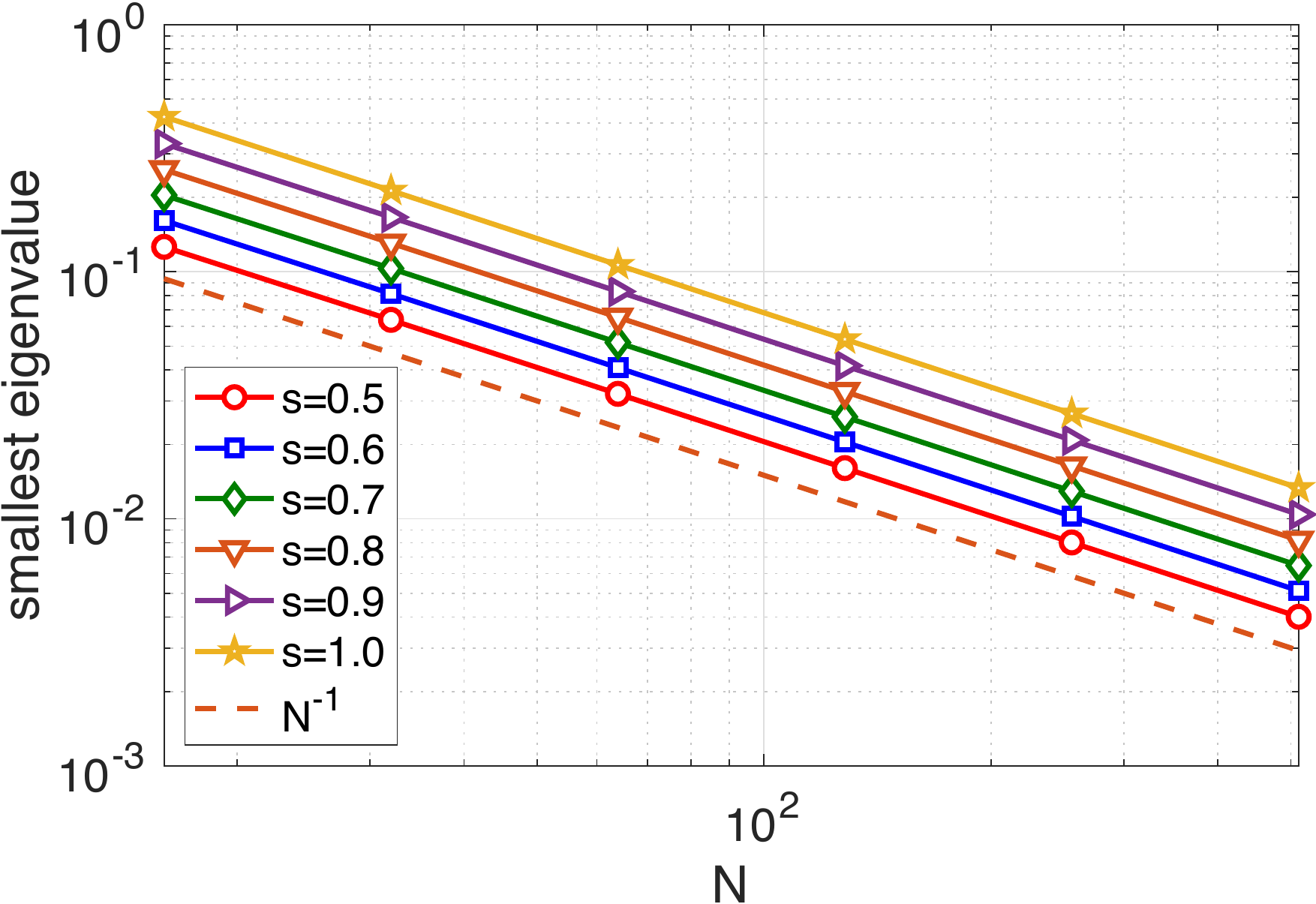}}
\end{minipage}}
\subfigure[${\rm Cond}\sim N^{(\alpha-1)(2s-1)\mu(s)+2s}$\hspace*{12pt}]{\hspace{-28pt}
\begin{minipage}[t]{0.40\textwidth}
\centering 
\rotatebox[origin=cc]{-0}{\includegraphics[width=0.8\textwidth,height=0.64\textwidth]{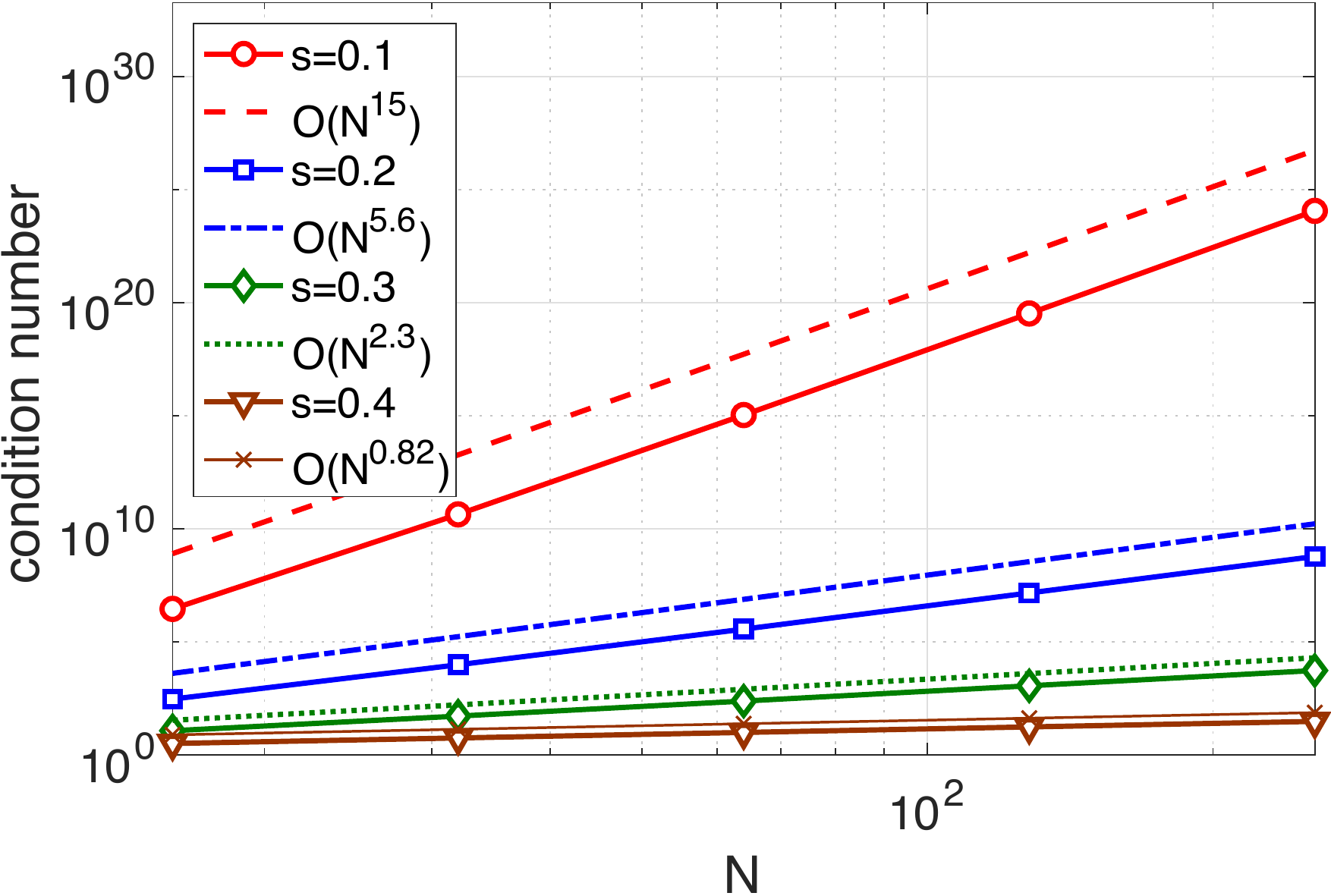}}
\end{minipage}}\hspace*{-35pt} {~}

\vskip -7pt
\caption
{\small Conditioning and  the smallest eigenvalue of the stiffness matrix $\bs S$ with $\alpha=2/s.$ 
(a)-(b):  various $s\in [1/2, 1].$ 
(c): various $s\in (0, 1/2)$.}\label{conds}
\end{figure}

The result  \eqref{condm}  is unknown for the fractional case.  Here, we explore this numerically, and provide some predictions or conjectures subject to rigorous proofs in future works. Note that in some critical situations (e.g.,  small $s$ or  very large $N$),   
we resort to the Multiprecision Computing Toolbox for Matlab \cite{MCT}.
 We highlight below the main numerical findings for $\bs S$ on the graded mesh generated by \eqref{xgrids}  with $\alpha=\beta>1,$
 and mostly consider $\alpha=2/s$ (the optimal value to achieve the best second-order accuracy for functions with algebraic endpoint singularities). 
 \begin{itemize}
 \item[(i)] The result \eqref{condm} is extendable to   $s\in [1/2, 1],$ that is, 
 \begin{equation}\label{condnum0}
{\rm Cond}(\bs S)=
 c\big(h_{\max}/h_{\min}\big)^{2s-1}N^{2s},  \;\; \,\, s\in[1/2,1]. 
\end{equation}
In Figure \ref{conds}(a), we illustrate the growth of the condition numbers with various $s\in[1/2,1],$   
and  find a good agreement between the numerical results and  \eqref{condnum0}. 
It is known that the smallest eigenvalue of $\bs S$  for the usual Laplacian (i.e.,  $s=1$) on a uniform mesh  behaves like 
${\lambda}_{\min}\approx\pi^2 h$ with $h$ being the mesh size.   Indeed, we observe from Figure \ref{conds}(b) that 
\begin{equation}\label{mineig}
{\rm \lambda}_{\min}(\bs S) =c h_{\max}=cN^{-1},  \, \,\,\,\, s\in[1/2,1]. 
\end{equation}
In fact, we also observe similar behaviours in \eqref{condnum0}-\eqref{mineig} for
various $\alpha>1,$ though we do not report the results here.  
\medskip 
\item[(ii)] The result \eqref{condnum0} does not hold for $s\in (0,1/2).$ We conjecture from numerical tests that 
\begin{equation}\label{condnum1}
{\rm Cond}(\bs S)=c\big(h_{\max}/h_{\min}\big)^{\mu(s)(1-2s)}\!N^{2s},  \, \,\,\,\, s\in (0, 1/2),\hspace{-8pt} 
\end{equation}
where $\mu(s)$ is some function. We refer to  Table \ref{TabN0} for some samples, and  find from ample tests that $\mu\in (0,1).$
We also observe from Figure\,\ref{conds}(c) that the condition number  increases rapidly as $s$ becomes smaller and closer to $0.$
\vspace*{-10pt}
\begin{table}[!th]
\centering\small
\caption{Samples of $\mu(s)$ with $\alpha=2/s$.} \label{TabN0}
\vspace*{-8pt}
\begin{tabular}{c c c c c c c c c}
\hline
\cline{1-1}
$s$           & $0.1$           & $0.15$        & $0.2$      & $0.25$      & $0.3$         & $0.35$       & $0.4$       & $0.45$ \\ \hline
$\mu(s)$   & $0.9505$    & $0.9394$    & $0.9032$ & $0.8369$ & $0.7167$    & $0.4868$   & $0.0113$ & $0.0144$ \\ \hline
\end{tabular}
\end{table}
 \end{itemize}

\vspace{-10pt}
 \subsection{Numerical results} 
\begin{figure}[!th]
\subfigure[{\rm Conv} $\sim N^{-\min\{2,\alpha s\}}$ \hspace*{12pt}]{\hspace{-20pt} 
\begin{minipage}[t]{0.40\textwidth}
\centering 
\rotatebox[origin=cc]{-0}{\includegraphics[width=0.8\textwidth,height=0.65\textwidth]{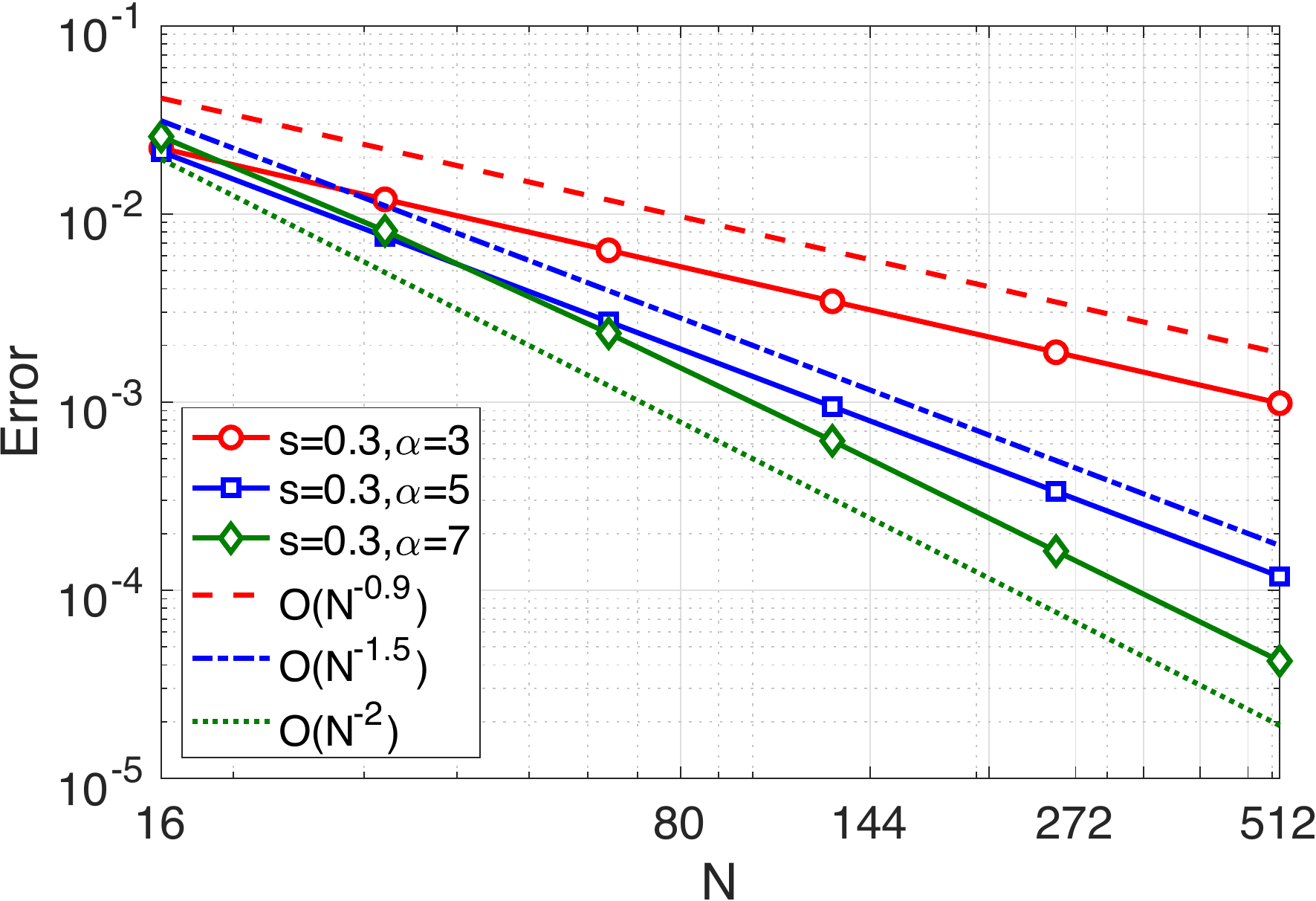}}
\end{minipage}}\hspace{-35pt}
\subfigure[{\rm Conv} $\sim N^{-2}, \alpha=2/s$  \hspace*{-20pt}]{ 
\begin{minipage}[t]{0.40\textwidth}
\centering
\rotatebox[origin=cc]{-0}{\includegraphics[width=0.8\textwidth,height=0.635\textwidth]{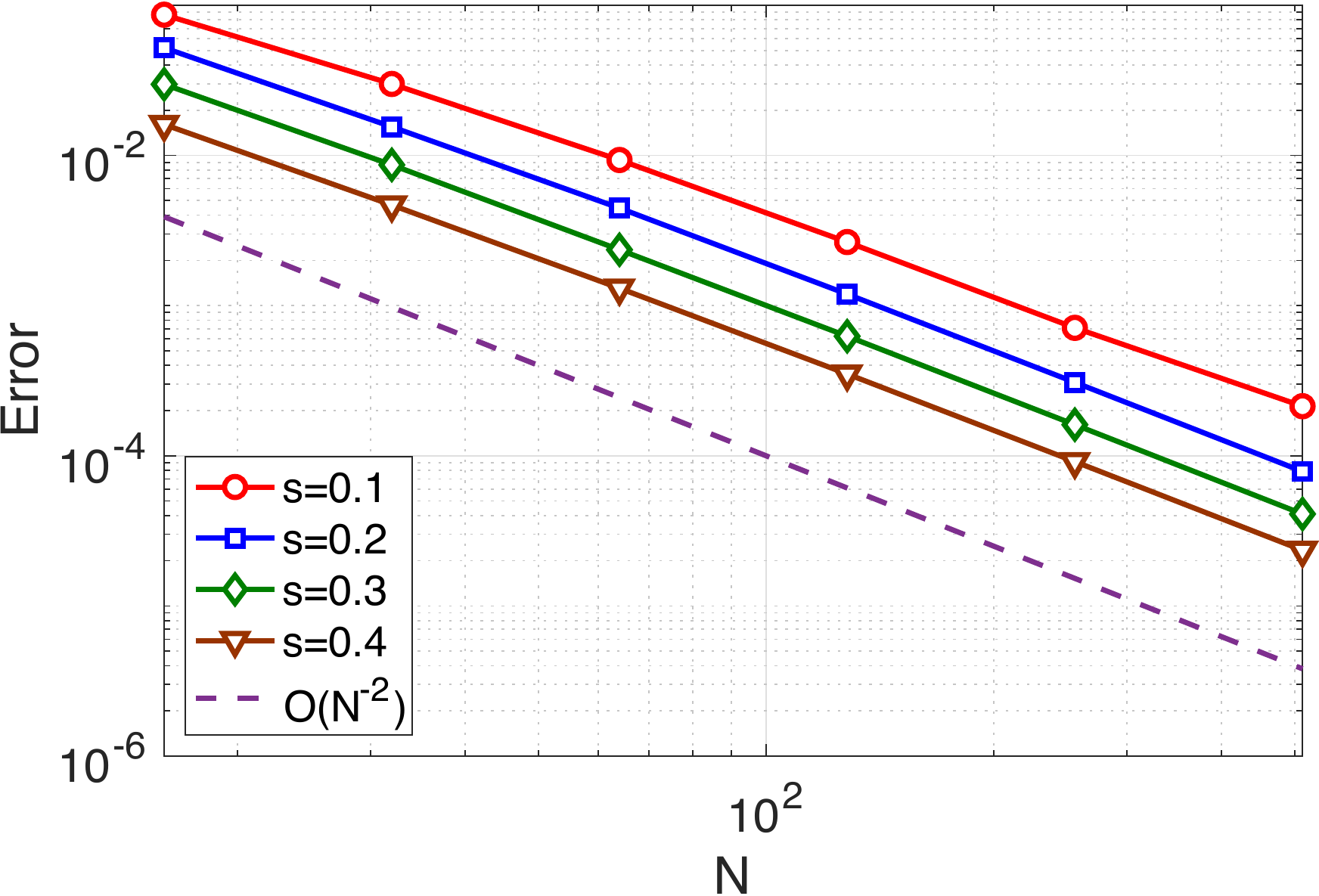}}
\end{minipage}}\hspace{-28pt}
\subfigure[{\rm Conv} $\sim N^{-2}, \alpha=2/s$ \hspace*{-30pt}]{\hspace{-18pt} 
\begin{minipage}[t]{0.40\textwidth}
\centering 
\rotatebox[origin=cc]{-0}{\includegraphics[width=0.8\textwidth,height=0.65\textwidth]{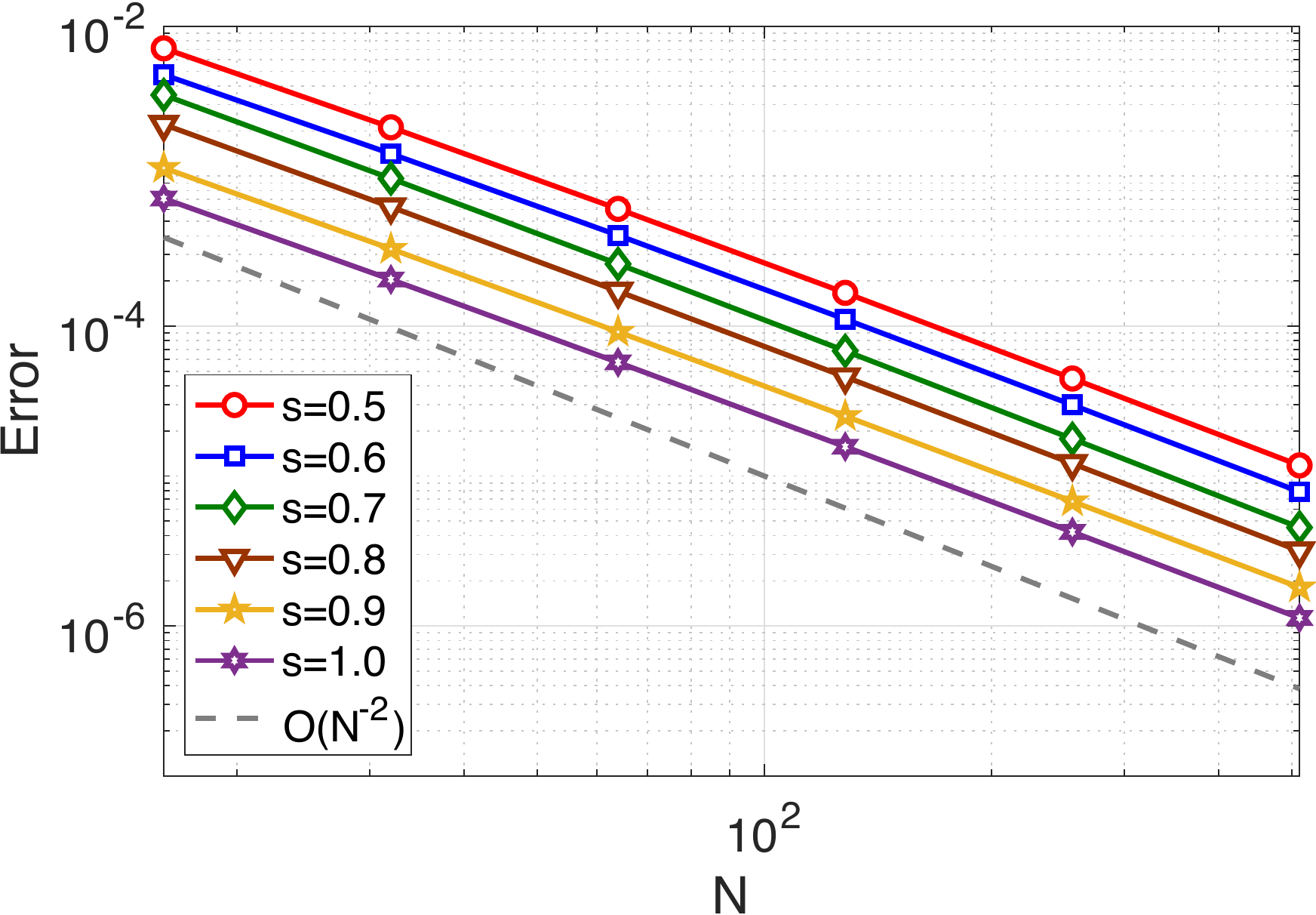}}\hspace*{-30pt} {~}
\end{minipage}}
\vskip -5pt
\caption
{\small Convergence order of the FEM solver on graded meshes.  (a): $s=0.3$ and different $\alpha$.  (b):  various  $s\in (0,1/2)$ and $\alpha=2/s$.  (c):    various  $s\in [1/2,1]$ and $\alpha=2/s$.}\label{err1}
\end{figure}

It is known that the solution of the fractional Poisson problem \eqref{dDsLap00} with a smooth source term $f(x)$ exhibits singularities near the boundary of $\Omega$ (cf.\,\cite{grubb2015fractional}).  In particular,   we find
 from \cite{dyda17fractional} that 
 \begin{equation*} 
 (-\Delta)^s ((1-x^2)^s_+ P_n^{(s,s)}(x))=\frac{\Gamma(n+2s+1)}{n!} P_n^{(s,s)}(x),\quad x\in \Omega=(-1,1), \;\; s>0, 
 \end{equation*}
 where $P_{n}^{(s,s)}(x)$ is the Jacobi polynomial of degree $n$, and  $u_+(x)=\max\{u(x), 0\}$. 
In the following computation, we take $f(x)=1$ in \eqref{dDsLap00}, and its exact solution is $ u(x)=(1-x^{2})^s_+/{\Gamma (2s+1)}$.
Following the same lines as in the proof of \cite[Thm.\,6.2.4]{MR2128285}, we can show that 
\begin{equation}\label{2rate0}
\max_{|x|\in \bar \Omega}\big|(u-I_{h}u)(x)\big| \le c N^{-\min\{2,\alpha s\}},
\end{equation}
where $I_{h} u$ is the piecewise linear FEM interpolation  of $u$ on the mesh  \eqref{xgrids} with $\alpha=\beta>1.$ As a result, 
the optimal order can be achieved when $\alpha=2/s.$  
With the explicit form of $\bs S$ in Theorem \ref{thmmain} and the aid of the Matlab toolbox \cite{MCT} 
(for some extreme situations, e.g., $s=0.1$ with ${\rm Cond}(\bs S)\sim N^{15},$ see Figure \ref{conds}(c)),  we demonstrate that 
the  same accuracy can be attained when  the FEM solution $u_h$  of \eqref{dDsLap02} is  in place of  $I_hu$ in \eqref{2rate0}. 
We observe from the numerical error plots in Figure \ref{err1}   that the convergence rate  of the FEM solver agrees well with the theoretical prediction.

\vspace{-15pt}
 \section{Concluding remarks and discussions}

 Different from the implementation of FEM  in the physical space, we computed the stiffness matrix of piecewise linear FEM for the IFL in the frequency space, and derived  the exact form of the entries. In fact, this approach can be extended to two-dimensional rectangular elements, but it is much more involved, which we shall report in a separate work.  Here, we  studied the graded mesh, and numerically demonstrated how the condition number of the stiffness matrix grew with the  parameters.  One  message is that computation with multiple precision is necessary, in order  to reduce the  round-off errors in  evaluating the entries of the fractional stiffness matrix and battling its large condition number.  
 In this study, we only considered the fractional order  $s\in (0,1],$ but the formulas in Theorem \ref{thmmain}  are valid for $s<3/2,$
which we leave for future investigation.

%
%
%
%

%
  
  
 \bibliographystyle{siam}

  \end{document}